# THE EXIT PROBLEM FOR DIFFUSIONS WITH TIME-PERIODIC DRIFT AND STOCHASTIC RESONANCE

By Samuel Herrmann and Peter Imkeller

*Université Henri Poincaré Nancy I and Humboldt-Universität zu Berlin*

Physical notions of stochastic resonance for potential diffusions in periodically changing double-well potentials such as the spectral power amplification have proved to be defective. They are not robust for the passage to their effective dynamics: continuous-time finite-state Markov chains describing the rough features of transitions between different domains of attraction of metastable points. In the framework of one-dimensional diffusions moving in periodically changing double-well potentials we design a new notion of stochastic resonance which refines Freidlin's concept of quasi-periodic motion. It is based on exact exponential rates for the transition probabilities between the domains of attraction which are robust with respect to the reduced Markov chains. The quality of periodic tuning is measured by the probability for transition during fixed time windows depending on a time scale parameter. Maximizing it in this parameter produces the stochastic resonance points.

**0. Introduction.** One of the simplest and earliest stochastic climate models goes back to [1] and [17]. It intends to give a qualitative explanation of glacial cycles and is based on a deterministic differential equation for the global mean temperature expressed through a balance between the albedo-driven absorbed and the black-body type emitted radiative energies. A periodic exterior forcing comes from the slowly fluctuating solar constant and is due to Milankovich cycles caused by the gravitation of big planets. Only the addition of a stochastic term as a second forcing makes spontaneous transitions between the otherwise isolated *metastable states* of temperature possible. The resulting stochastically and periodically perturbed differential equation was capable of describing at least one characteristic aspect of experience: the typically short and abrupt *transitions*, observed before by Kramers [14] in reaction-diffusion phenomena. The model was soon strongly disputed. Despite its lack of realistic assumptions, the concepts underlying the model brought to light the phenomenon of *stochastic resonance*. Roughly speaking, a periodic (input) system subject to random perturbations is in stochastic resonance, if the noise intensity is tuned in such a way that the random periodic output is optimal. A very lively research field developed around this concept, drawing numerous examples from a wide spectrum of areas (see [9, 10] for a survey).

The mathematically precise understanding of the phenomenon is still under discussion at the time this paper is written. The first approach is done in [7], where the deep large deviations' theory of [8] is employed to produce a notion of stochastic resonance explaining the phenomenon in the small noise limit as approximating the periodic hopping between the energetically most favorable states in the landscape provided by a periodically weakly perturbed potential with finitely many local minima. In this sense stochastic resonance can be understood as the ability of the system to undergo quasi-periodic motion in the limit of small noise intensity. Let us briefly recall this interpretation more precisely. If noise intensity is $\varepsilon$, in the absence of periodic exterior forcing, the exponential order of times at which successive transitions between metasable states occur corresponds to the work to be done against the potential gradient to leave a well. This fact, heuristically derived by Kramers and Eyring (*Kramers' time*), is shown with mathematical accuracy in [8]. The attractor basins are subdivided into a hierarchy of cycles with main states corresponding to the deepest among the cycle states. In the presence of periodic forcing with period time scale $e^{\mu/\varepsilon}$, in the limit $\varepsilon \to 0$ transitions between (the main states of) cycles with critical hopping work close to $\mu$ will be periodically observed. Transitions with smaller critical work may happen, but are negligible in the limit. Those with larger critical work are forbidden. In the simplest case of two minima of potential depth $V$ and $v$, $v < V$, the role of which switches periodically at time $T$, for $T$ larger than $e^{v/\varepsilon}$ the diffusion will be quasi-deterministic, that is, close to the deterministic periodic function following the location of the deepest well.

Quasi-periodicity captures an important aspect of stochastic resonance, as it provides conditions under which stochastic trajectories are able to exhibit periodic behavior. Yet, physics literature (see [9, 10]) stipulates that stochastic resonance not only explains conditions for *stochastically periodic behavior* but comprises its optimality in a sense quite similar to the resonance notions of wave dynamics. In classical optics resonance is understood as the optimal amplitude of the response of the system to periodic excitation. In the same way, a stochastic resonance point is claimed to explain *optimal periodic tuning* of the stochastic trajectories of the diffusion responding to deterministic periodic excitation. Amplitude as a measure of quality of periodic tuning is replaced by *signal-to-noise ratio* or *spectral power amplification* (see below). Numerical simulations as, for example, in [16] clearly support the optical evidence that beyond the threshold described by Freidlin [7] at which quasi-deterministic behavior becomes possible, for different noise intensities quite different qualities of periodicity of the random trajectories can be observed. There are parameter ranges for $T$ in which random trajectories follow quite well the deterministic shapes of excitation curves. But as $T$ gets even bigger, many short excursions to the wrong well during one period may occur. They will not count on the exponential scale on which quasi-periodic motion is measured, but trajectories will look less and less periodic. Physicists' quality measures for tuning therefore cannot be explained on the basis of quasi-deterministic motion alone.

The thesis by Pavlyukevich [13] and [18] presents an attempt to provide a mathematically sound underpinning of physical notions of stochastic resonance based on optimality of periodic tuning—as opposed to the trajectorial analysis of Berglund and Gentz [2] containing very fine estimates on relaxation times. The physical concepts are mostly based on comparisons of trajectories of the noisy system and the deterministic periodic curve describing the location of the relevant metasable states, averaged with respect to the equilibrium measure. In the simple one-dimensional situation considered above the system switches between a double-well potential state $U(x)$ with two wells of unequal depths $V$ and $v$, $v < V$, during the first half period, and $U(-x)$ for the second half period. The total period length is $T$, and stochastic perturbation comes from the coupling to a white noise of intensity $\varepsilon$. The most important measures of quality studied are the *spectral power amplification*, the related *signal-to-noise ratio* or the *entropy of the equilibrium distribution*. In particular, the first two mentioned play an eminent role in the physical literature. They mainly contain the $L^2$ average in equilibrium of the spectral component of the solution trajectories corresponding to the input period $T$, normalized in different ways. These measures of quality are functions of $\varepsilon$ and $T$, and the problem of finding the resonance point, for example, consists in optimizing them in $\varepsilon$ for fixed (large) $T$.

Let us briefly explain a striking shortcoming of these concepts of optimal periodic tuning which made us look for different ones. The first step to find optimal tuning intensities $\varepsilon(T)$ for large $T$ consists in reducing the dynamics of the diffusion to the *interwell motion*, that is, the pure transitions between the potential minima. In the physics literature, this corresponds to the reduction given, for example, by [15]. One ends up with continuous-time two-state Markov chains with transition probabilities corresponding to the inverses of the diffusions' Kramers–Eyring times. The mathematical analysis of stochastic resonance then proceeds along the following lines. One first determines the optimal tuning parameters $\varepsilon(T)$ for large $T$ for the approximating Markov chains, a rather simple task. To see that the Markov chain's behavior approaches the diffusion's in the small noise limit, spectral theory of the infinitesimal generator is used. Its spatial part is seen to possess a spectral gap between the second and third eigenvalues, and therefore the closeness of equilibrium distributions of the Markov chain on the one hand and diffusion on the other hand can be well controlled. Surprisingly, however, the notion of spectral power amplification is not robust for the passage from the Markov chain to the diffusion. Subtle dependencies on the geometrical fine structure of the potential function $U$ in the minima beyond the expected curvature properties lead to quite unexpected counterintuitive effects. For example, a subtle drag away from the other well caused by the sign of the third derivative of $U$ in the deep well suffices to make the spectral power amplification curve strictly increasing in the parameter range in which the approximating Markov chain has its resonance point. This dramatic deviation from expected behavior is due to the significance the spectral power amplification attributes to small *intrawell fluctuations*.

Our main motivation in writing this paper was to investigate concepts by which on the one hand the physical intuition of optimal periodic tuning of random trajectories with a simple periodic input can be interpreted in a mathematically sound way, and which on the other hand fail to have this unfortunate defect of robustness. We deal with the framework of one-dimensional potential diffusions. The notion of quality of periodic tuning we shall investigate completely excludes the effect of small intrawell fluctuations and purely relies on the transition mechanism between domains of attraction given by the potential. At the same time it generalizes the previously known results to potential functions which may vary periodically in time in a continuous, but otherwise quite general way, and whose growth at $\pm\infty$ may just be linear. More precisely, we study diffusion processes driven by a Brownian motion of intensity $\varepsilon$ given by the stochastic differential equation (SDE)

$$dX_t = -\frac{\partial}{\partial x} U\left(\frac{t}{T}, X_t\right) dt + \sqrt{2\varepsilon}\, dW_t, \qquad t \geq 0.$$

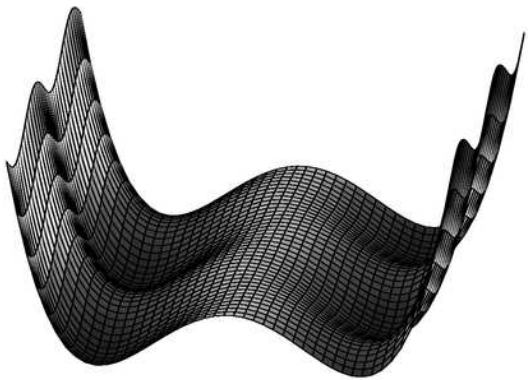

Fig. 1. *Potential landscape.*

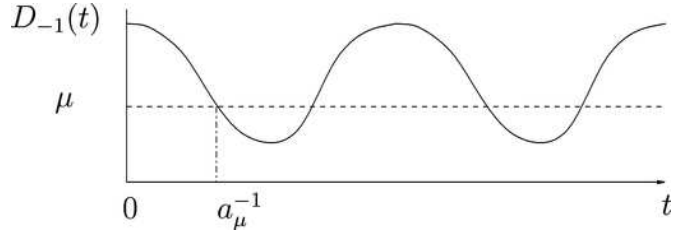

Fig. 2. *Definition of $a_\mu^{-1}$.*

The underlying potential landscape (see Figure 1) is described by a function $U(t,x), t \geq 0, x \in \mathbb{R}$, which is periodic in time with period 1, and its temporal variation, by the rescaling with very large $T$, acts on the diffusion at a very small frequency. $U$ is supposed to have exactly two wells located at $\pm 1$, separated by a saddle at 0. The depth (measured in positive quantities) of $U(t, \cdot)$ at $\pm 1$ is given by the 1-periodic depth functions $D_{\pm 1}(t)$ which are assumed to never fall below zero. We shall throughout look at time scales $T = \exp(\frac{\mu}{\varepsilon})$, for which the Kramers–Eyring formula indicates that transitions, for example, from the domain of attraction of $-1$ to the domain of attraction of 1 will occur as soon as $D_{-1}$ becomes less than $\mu$, that is, at time (see Figure 2)

$$a_\mu^{\pm 1} = \inf\{t \geq 0 : D_{\pm 1}(t) \leq \mu\}.$$

This triggers periodic behavior of the diffusion trajectories on long time scales. The modern theory of metastability in *time-homogeneous diffusion processes* complements the fundamental large deviations' theory presented by [8] to produce the exponential decay rates of transition probabilities between different domains of attraction of a potential landscape together with very sharp multiplicative error estimates, uniformly on compacts in system parameters. Their sharpest forms are presented in some papers by Bovier, Eckhoff, Gayrard and Klein [3, 4], improving Day's previous results obtained in [5, 6]. They are derived from deep relationships of large deviations' theory with the spectral and capacity theory of the infinitesimal generator. We shall make use of this powerful machinery to obtain very precise estimates of the exponential tails of the laws of the transition times between domains of attraction. In fact, we have to extend the estimates by Bovier, Gayrard and Klein [4] to the framework of *time-inhomogeneous diffusions*. In the underlying one-dimensional situation, this is roughly achieved by freezing the time dependence of the potential on small time intervals to define lower and upper bound time-homogeneous potentials not differing very much from the original one. Consequently, comparison theorems are used to control the transition behavior from above and below by the behavior of the corresponding time-homogeneous diffusions. This allows very precise estimates on the probabilities with which the diffusion at time scale $T = \exp(\frac{\mu}{\varepsilon})$ transits from the domain of attraction of $-1$ to the domain of attraction of $1$ or vice versa within time windows $[(a_\mu^{\pm 1} - h)T, (a_\mu^{\pm 1} + h)T]$ for small $h > 0$. If $\tau_x(X)$ denotes the transit time to $x$, it is shown to be given by

$$\lim_{\varepsilon \to 0} \varepsilon \ln(1 - M(\varepsilon, \mu)) = \max_{i=\pm 1}\{\mu - D_i(a_\mu^i - h)\},$$

with

$$M(\varepsilon, \mu) = \min_{i=\pm 1} \mathbb{P}_i(\tau_{-i}(X) \in [(a_\mu^i - h)T, (a_\mu^i + h)T]), \qquad \varepsilon > 0, \mu \in I_R,$$

and where $I_R$ is the *resonance interval* (Figure 3), that is, the set of scale parameters for which trivial or chaotic transition behavior of the trajectories is excluded (Figure 4).

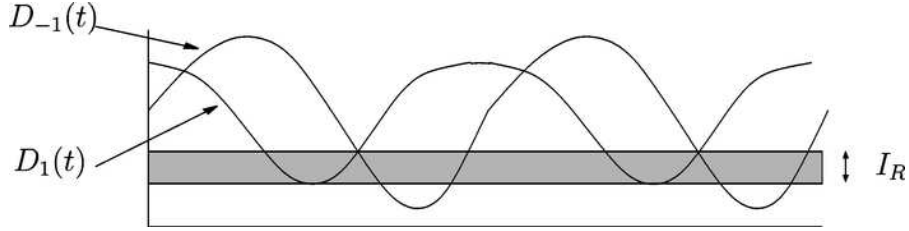

Fig. 3. *Resonance interval.*

The stated convergence is *uniform* in $\mu$ on compact subsets of $I_R$. This allows us to take $M(\varepsilon, \mu)$ as our measure of periodic tuning, compute the scale $\mu_0(h)$ for which the transition rate is optimal, and define the *stochastic resonance point* as the eventually existing limit of $\mu_0(h)$ as $h \to 0$. This notion of stochastic resonance is strongly related to the notions of periodic tuning based on *interspike intervals* (see Figure 5 and [11]), which describe the probability distribution for transitions as functions of time with exponentially decaying *spikes* near the integer multiples of the forcing periods. As opposed to the physics notions based on spectral decomposition of the statistics of the solution trajectories investigated in [13] or [18] it has the big advantage of being robust for the passage from the diffusion to the two-state Markov chain reducing its behavior to the features of pure transitions between the two domains of attractions of metasable points.

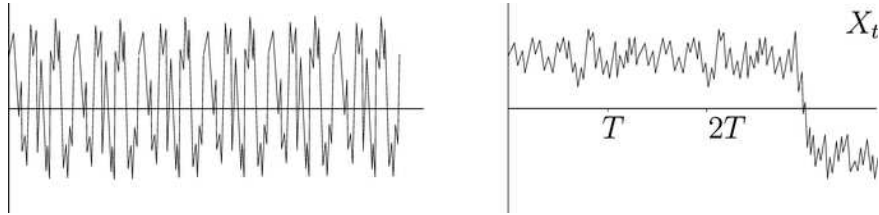

Fig. 4. *Chaotic and trivial transition behavior of the trajectories.*

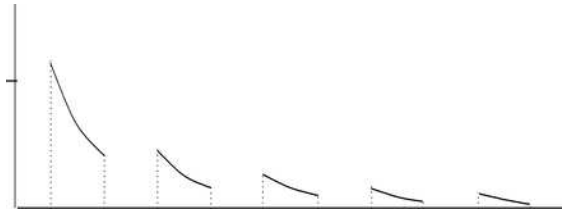

Fig. 5. *Length distribution of the interspike intervals in a simplified model (two-state Markov chain).*

Here is an outline of the organization of the material in the paper. Section 1 presents a review of results from the asymptotic theory of time-homogeneous diffusions and their metasable sets needed for our purposes (Theorem 1.1). In Section 2 we bring to work the tools of comparison theorems to deduce the sharp exponential transition rates for our time-periodic diffusions from the time-homogeneous results (Theorem 2.1). Section 3 is devoted to applying these sharp estimates to identify stochastic resonance points for diffusions (Theorem 3.2), compare them to their counterparts for the reduced Markov chains and prove robustness of our notion of resonance (Theorem 3.4).

**1. Exponential distribution of transition times for time-homogeneous diffusions.** It will turn out to be crucial for our approach of periodic tuning to be discussed later to obtain large deviation type estimates for the exponential decay rate of the law of transition times uniformly in a time scale parameter. We shall make use of a technique of freezing time-dependent potentials on small subintervals of the periodicity interval $[0, T]$ on their states taken at fixed times in the intervals, to be able to use known results for time-homogeneous diffusions. In this setting, the uniformity problem translates into uniformity of the convergence to exponential decay rates in compact subsets of the domain of attraction the diffusion starts in and in time. It is clear that we are led directly into large deviations' estimates for exit time distributions of *time-homogeneous diffusions* such as presented in the pioneering book by Freidlin and Wentzell [8]. But for obtaining uniformity in space and time, one has to use sharpened versions of these estimates developed later for controlling in particular the exponential errors in the estimates. The purpose of this section is to summarize what we shall need from this fine well-established theory.

We shall refer to the most recent and advanced development of sharp estimates for transition times presented in [3, 4]. They are valid far beyond our modest framework, both in the multidimensional case and for any finite number of local minima of the potential. Their quality comes from a detailed analysis of the relationship between transition times and low-lying eigenvalues of the spectrum of the infinitesimal generator of the diffusion. We shall state them in the simple one-dimensional setting given here. A more complex multidimensional version can also be found in [5]. For this purpose, suppose that $Q$ is a purely space-dependent $\mathcal{C}^2$ potential function (see Figure 6) possessing only $-1, 1$ as local minima, separated by the saddle point 0 at which $Q$ takes the value 0. Suppose that the curvature of $Q$ at $-1$ is strictly positive, that is, $Q''(-1) > 0$. As for ultra- or hypercontractivity type properties for $Q$, we shall assume that it has exponentially tight level sets; that is, there is $a_0 > 0$ such that for any $a \geq a_0$ there exists a constant $C(a)$ such that for $\varepsilon \leq 1$,

$$(1) \qquad \int_{\{y:Q(y)\geq a\}} \exp\left(-\frac{Q(z)}{\varepsilon}\right) dz < C(a) \exp\left(-\frac{a}{\varepsilon}\right).$$

We shall concentrate in this situation on a transition out of the domain of attraction of the stable point $-1$ for the diffusion associated with the SDE

$$dY_t^\varepsilon = -Q'(Y_t^\varepsilon)\, dt + \sqrt{2\varepsilon}\, dW_t,$$
$$Y_0^\varepsilon = y.$$

Let $C$ be a closed interval of the form $[d, \infty[$ with $d \neq 0$. To state our aim in a slightly different version, we will be interested in the asymptotics of the entrance time of $Y^\varepsilon$ into $C$:

$$\tau_C^\varepsilon = \inf\{t > 0 : Y_t^\varepsilon \in C\}.$$

Then we obtain the following result (see [4] or [5]).

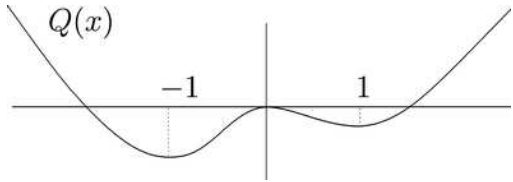

Fig. 6. *Potential.*

THEOREM 1.1. *Let $\lambda^\varepsilon$ denote the principal eigenvalue of the linear operator*

$$\mathcal{L}^\varepsilon u = \varepsilon u'' - Q' u'$$

*with Dirichlet boundary conditions on $\partial C$. Then for every compact $K \subseteq {]{-\infty, 0}[}$ there is a constant $c > 0$ such that*

(2) $$\mathbb{P}_y(\tau_C^\varepsilon > t) = e^{-\lambda^\varepsilon t}(1 + \mathcal{O}_K(e^{-c/\varepsilon})),$$

*where $\mathcal{O}_K$ denotes an error term which is uniform in $y \in K$, $t \geq 0$. Moreover, for the asymptotic behavior of the eigenvalue $\lambda^\varepsilon$ the following holds*:

(3) $$\lambda^\varepsilon \mathbb{E}_y[\tau_C^\varepsilon] \to 1 \qquad \textit{uniformly on compacts } K \subseteq {]{-\infty, 0}[}$$

*as $\varepsilon \to 0$.*

PROOF. There are two small issues which deserve comments.

First, the uniformity over compacts in $]{-\infty, 0}[$ claimed in the main statements. Day [5] tackles it. But he considers only exits from bounded domains. Bovier, Gayrard and Klein [4] have a version for unbounded domains, but uniformity over compacts of the domain of attraction in which the diffusion starts is not explicitly proved. It is, however, hidden in their method of proof of Theorem 1.3 ([4], pages 30 and 31) which makes use of an eigenfunction expansion. But due to regularity results on the eigenfunctions (see [4], pages 16–18) they must be bounded on compacts in the domain of attraction $]{-\infty, 0}[$. This implies the desired uniformity.

The second comment concerns our assumptions on $C$. Translated into our setting, in [4] the target set $C$ is assumed to be closed, to contain a neighborhood of 1 if the potential is deeper there than at $-1$, and to have a positive distance from the saddle 0. Since we are in a one-dimensional setting, we may reduce these conditions to the simple one $d \neq 0$. If necessary, we may always cut out of $C$ a small open neighborhood of 0 without changing the law of $\tau_C$ if starting from the domain of attraction of $-1$. $\square$

**2. Exponential transition rates between moving domains of attraction.** We shall now consider a potential diffusion given by the one-dimensional SDE

(4) $$dX_t = -\frac{\partial U}{\partial x}\left(\frac{t}{T}, X_t\right) dt + \sqrt{2\varepsilon}\, dW_t.$$

The time-periodic potential $U$ of period 1 is supposed to fulfill the following conditions. First of all, its global rough geometry is the one of a double-well potential with temporally moving wells, but time-independent critical points. For simplicity we suppose that its local minima are given by $\pm 1$, and its only saddle point by 0, independently of time. So $\pm 1$ are the only metastable states of the potential on the whole time axis. Outside of $0, \pm 1$, $\frac{\partial U}{\partial x}$ is supposed to be continuous in $(t, x)$. Our main concern will be the asymptotics of the transition times from the domain of attraction $]-\infty, 0[$ of $-1$ to the domain of attraction $]0, \infty[$ associated with 1 of the *time-inhomogeneous diffusion* in the small noise limit $\varepsilon \to 0$. More precisely, we will be interested in describing the exponential transition rate from the domain of $-1$ to the domain of 1. Our potential not being time-homogeneous, we shall make use of comparison arguments with diffusions possessing time-independent potentials in order to find a careful reduction of the inhomogeneous exit problem to the homogeneous one, and use the asymptotic results stated in Theorem 1.1 in this framework. This will be achieved by freezing the driving force derived from the potential on small time intervals on the mimimal or maximal level it takes there. To be more precise, for each interval $I \subset \mathbb{R}_+$ let

$$(5) \quad V_I(x) = \sup_{t \in I} \frac{\partial U}{\partial x}(t, x) \quad \text{and} \quad R_I(x) = \inf_{t \in I} \frac{\partial U}{\partial x}(t, x).$$

See Figure 7. The regularity conditions valid for $U$ imply that $V$ and $R$ are continuous functions. Moreover, $V_I(-1) = R_I(-1) = 0$. If $I = [a, b]$, we denote by $\overline{X}^I$ the solution of the SDE on $\mathbb{R}_+$

$$(6) \quad \begin{aligned} d\overline{X}_t^I &= -R_I(\overline{X}_t^I)\, dt + \sqrt{2\varepsilon}\, dB_t, \\ \overline{X}_0^I &= X_{aT}. \end{aligned}$$

$\underline{X}^I$ is defined in the same way, replacing $R_I$ by $V_I$. These two *time-homogeneous* diffusions are used to control the *time-inhomogeneous* diffusion $X$ as long as time runs in the interval $I$. In fact, we have, $P$-a.s.,

$$\underline{X}_{tT}^I \leq X_{(t+a)T} \leq \overline{X}_{tT}^I, \qquad t \in [0, b-a].$$

Of course, to make use of the asymptotic results stated in the previous section, we need ultra- or hypercontractivity properties for the frozen potentials. To formulate a hypothesis which is both not too restrictive and easy to handle for time-dependent potentials, let us give the following easy sufficient criterion for exponential tightness of levels of a time independent potential $Q$.

LEMMA 2.1. *Assume that $Q$ is a real-valued differentiable function on $\mathbb{R}$, and that there are constants $K_1, K_2 > 0$ such that*

$$(7) \quad \begin{aligned} Q'(x) &\leq -K_2 \quad \text{for } x \leq -K_1, \\ Q'(x) &\geq K_2 \quad \text{for } x \geq K_1. \end{aligned}$$

*Then $Q$ has exponentially tight level sets.*

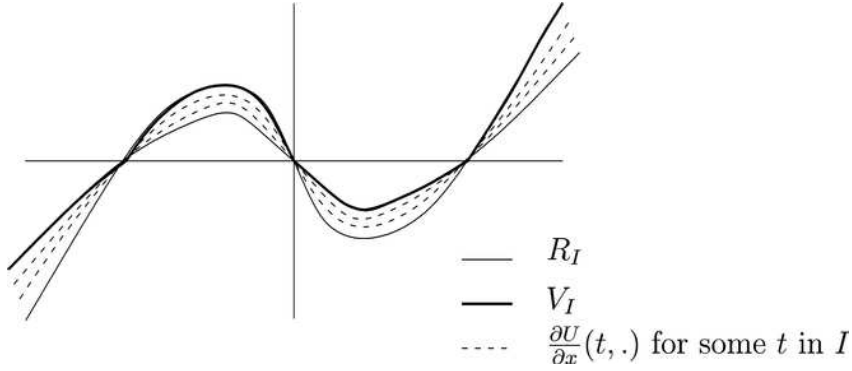

Fig. 7. *Definition of $V_I$ and $R_I$.*

PROOF. It is obviously enough to argue on $\mathbb{R}_+$. Due to (7), we know that near $\infty$, $Q$ is strictly increasing with inverse $Q^{-1}$. So for $a > 0$ big enough and $\varepsilon > 0$ we have

$$\int_{\{y:Q(y)\geq a\}} \exp\left(-\frac{Q(z)}{\varepsilon}\right) dz = \int_a^\infty \exp\left(-\frac{t}{\varepsilon}\right) \frac{1}{Q'(Q^{-1}(t))} dt$$
$$\leq \frac{1}{K_2} \int_a^\infty \exp\left(-\frac{t}{\varepsilon}\right) dt$$
$$= \frac{\varepsilon}{K_2} \exp\left(-\frac{a}{\varepsilon}\right).$$

This clearly implies exponential tightness. □

Conditions like (7) are practical in our setting for the following reason: if we assume them to be satisfied uniformly for $t \in [0,1]$, it is clear that all frozen potential functions will inherit the property. We will therefore assume throughout that there exist constants $K_1, K_2 > 0$ such that

(8)
$$\sup_{t \geq 0} \frac{\partial U}{\partial x}(t,x) \leq -K_2 \qquad \text{for } x \leq -K_1,$$
$$\inf_{t \geq 0} \frac{\partial U}{\partial x}(t,x) \geq K_2 \qquad \text{for } x \geq K_1.$$

We measure periods $T$ on the logarithmic scale $\mu$ given by $T = T_\varepsilon = \exp(\frac{\mu}{\varepsilon})$. The reason for this is hidden in the classical formula of Kramers–Eyring. It states that in the small noise limit $T$ is the time it takes the diffusing particle to climb a height $\mu$ in a potential landscape $U$. This formula has the following intuitive consequences. Assume the diffusion faces an obstacle of constant potential height $U^+ > \mu$ for exiting if it diffuses on time scale $T$. Then asymptotically it never exits on this scale. On the other hand, if it faces an obstacle of height $U_- < \mu$ diffusing at time scale $T$, it has to exit immediately in the small noise limit. We introduce the depth function at $x \in \mathbb{R}$ by

$$D_x(t) = U(t,0) - U(t,x), \qquad t \geq 0.$$

The maximal well depths $D_{-1}$ and $D_1$ will be of particular importance. We shall assume that they satisfy the assumption (see Figure 2)

(9)      all local extrema of $D_1$, $D_{-1}$ are global and the functions are strictly monotonous between the extrema.

For $\mu \geq 0$ let now

$$a_\mu = \inf\{t \geq 0 : D_{-1}(t) \leq \mu\}.$$

This is the same function as $a_\mu^{-1}$ defined above. Here we omit the superscript since we always concentrate on transitions from $-1$ to $1$. The interval

$$(10) \qquad I_{-1} = \left] \inf_{t \geq 0} D_{-1}(t), \sup_{t \geq 0} D_{-1}(t) \right[$$

contains all possible depths the potential minimum located at $-1$ takes during one period of time. Note that $a_\mu = 0$ for $\mu$ strictly above the upper boundary of $I_{-1}$, and $a_\mu = \infty$ for $\mu$ strictly below the lower boundary. If $\mu \in I_{-1}$, then for times beyond $a_\mu$, the barrier height between $-1$ and $0$ has dropped below the critical level so that on time scale $T$ the diffusion exits immediately. According to this heuristic, the diffusion running on time scale $T$ should exit the domain of attraction $]-\infty, 0[$ through $0$ and then transit to the other well immediately after $TD_{-1}(t)$ drops below $\mu$. We shall be interested in the exponential rate at which this happens, uniformly in starting points taken from a compact in $]-\infty, 0[$. For this purpose, for a regular diffusion $Y$ on $\mathbb{R}$ and $a \in \mathbb{R}$, we denote by $\tau_a(Y)$ the first hitting time of $a$ by $Y$.

Our first aim is to prove the precise estimate

$$(11) \qquad \lim_{\varepsilon \to 0} \varepsilon \ln \mathbb{P}_x(\tau_1(X) \leq (a_\mu - h)T) = \mu - D_{-1}(a_\mu - h)$$

for $x < 0$, uniformly in $\mu$ on compact subintervals $\Gamma$ of (10) and for $0 < h < \inf_{\mu \in \Gamma} a_\mu$. Note that given such a compact $\Gamma$, due to the continuity of $D_{-1}$ we have $\inf_{\mu \in \Gamma} a_\mu > 0$. Fix $\Gamma, x$ and $h$ from now on. Our approach proceeds in essentially two steps.

2.1. *Upper bound for the exponential exit rate.* In the first step we shall find upper bounds for the exponential rates. For this purpose we shall partition the relevant time intervals $[0, a_\mu - h]$. Fix some $\delta > 0$ such that $\delta < |x|$. Since $\frac{\partial U}{\partial x}$ is continuous, we may choose an equidistant partition $0 = r_0 < r_1 < \cdots < r_n = a_\mu - h$ of $[0, a_\mu - h]$ with mesh $\gamma$ small enough to ensure

$$(12) \qquad \sup_{s,t \in [r_j, r_{j+1}]} \sup_{x \in [-1,0]} \left| \frac{\partial U}{\partial x}(t,x) - \frac{\partial U}{\partial x}(s,x) \right| \leq \delta.$$

Denote $I_j = [r_{j-1}, r_j]$, $1 \leq j \leq n$. Though the choice of the intervals depends on $h$ and $\Gamma$, their number will be bounded by a universal constant. Using this partition, we may start our search for an upper bound by freezing the time-dependent potential on $I_j$ at its value $R_{I_j}$, $1 \leq j \leq n$, and then comparing the diffusion $X$ there with $\overline{X}^{I_j}$. There is a little difficulty with this procedure. The drift coefficients $R_{I_j}$ which govern the SDE for $\overline{X}^{I_j}$ were defined by taking infima over time intervals. This operation may destroy their differentiability properties in the spatial variable. Therefore it may be necessary to compare $X$ on the intervals $I_j$ with smoother diffusions still dominating it. But this can be done at no extra cost. For each $1 \leq j \leq n$ we may choose $R_j \in \mathcal{C}^1(\mathbb{R})$ satisfying

$$(13) \quad \begin{aligned} &R_j \leq R_{I_j} \leq R_j + \delta, \\ &\text{there are } \underline{m}_j \in\, ]-1-\delta, -1+\delta[\,,\ s_j \in\, ]-\delta, \delta[\,,\ \overline{m}_j \in\, ]1-\delta, 1+\delta[ \\ &\text{such that } R_j|_{]-\infty, \underline{m}_j[ \cup\, ]s_j, \overline{m}_j[} < 0,\ R_j|_{]\underline{m}_j, s_j[ \cup\, ]\overline{m}_j, \infty[} > 0, \\ &R'_j(\underline{m}_j) > 0, \\ &\text{there are constants } K_1, K_2 > 0 \text{ such that } R_j(x) \leq -K_2 \text{ for } x \leq -K_1, \\ &R_j(x) \geq K_2 \text{ for } x \geq K_1. \end{aligned}$$

Let $X^j$ be the diffusion associated with $R_j$, $1 \leq j \leq n$. Let us choose a partition $x = x_0 < x_1 < \cdots < x_n = -\delta$ of the interval $[x, -\delta]$ which will typically not be supposed to be equidistant. In the following inequality the diffusions $X$ and $X^j$ on $I_j$ are compared and the Markov property is employed. We have

$$\mathbb{P}_x(\tau_1(X) \leq (a_\mu - h)T)$$

$$\leq \sum_{j=1}^{n-1} \mathbb{P}_x(\tau_{x_{j-1}}(X) \geq r_{j-1}T, \tau_{x_j}(X) \leq r_j T)$$

$$+ \mathbb{P}_x(\tau_{x_{n-1}}(X) \geq r_{n-1}T, \tau_1(X) \leq r_n T)$$

$$\leq \sum_{j=1}^{n-1} \mathbb{P}_x[\{\tau_{x_j}(X_{\cdot + r_{j-1}T}) \leq \gamma T\} \cap \{\tau_{x_{j-1}}(X) \geq r_{j-1}T\}]$$

(14)
$$+ \mathbb{P}_x[\{\tau_1(X_{\cdot + r_{n-1}T}) \le \gamma T\} \cap \{\tau_{x_{n-1}}(X) \ge r_{n-1}T\}]$$

$$\le \sum_{j=1}^{n-1} \mathbb{E}_x[\mathbb{P}_{X_{r_{j-1}T}}(\tau_{x_n}(X^j) \le \gamma T)\mathbb{1}_{\{\tau_{x_{j-1}}(X) \ge r_{j-1}T\}}]$$

$$+ \mathbb{E}_x[\mathbb{P}_{X_{r_{n-1}T}}(\tau_1(X^n) \le \gamma T)\mathbb{1}_{\{\tau_{x_{n-1}}(X) \ge r_{n-1}T\}}]$$

$$\le \sum_{j=1}^{n-1} \mathbb{P}_{x_{j-1}}(\tau_{x_j}(X^j) \le \gamma T) + \mathbb{P}_{x_{n-1}}(\tau_1(X^n) \le \gamma T).$$

Let us now fix $1 \le j \le n$ and continue estimating the terms $\mathbb{P}_{x_{j-1}}(\tau_{x_j}(X^j) \le \gamma T)$ and $\mathbb{P}_{x_{n-1}}(\tau_1(X^n) \le \gamma T)$ individually. For this purpose we apply Theorem 1.1 for $Q = R_j, d = x_j$ to obtain that

(15)
$$\mathbb{P}_{x_{j-1}}(\tau_{x_j}(X^j) \le \gamma T)$$
$$\le 1 - e^{-\lambda_j^\varepsilon \gamma T}(1 - e^{-c/\varepsilon}) \le 1 - e^{-\lambda_j^\varepsilon \gamma T} + e^{-c/\varepsilon},$$
$$\mathbb{P}_{x_{n-1}}(\tau_1(X^n) \le \gamma T)$$
$$\le 1 - e^{-\lambda_n^\varepsilon \gamma T}(1 - e^{-c/\varepsilon}) \le 1 - e^{-\lambda_n^\varepsilon \gamma T} + e^{-c/\varepsilon},$$

uniformly in $T$, hence uniformly in $\mu$. Here $\lambda_j^\varepsilon$ denotes the principal eigenvalue of the operator $\mathcal{L}_j^\varepsilon$ defined by

$$\mathcal{L}_j^\varepsilon u = \varepsilon u'' - R_j u'$$

with Dirichlet boundary conditions at $x_j$. We now come to the crucial part of the derivation of an upper estimate. We shall use precise asymptotics of the eigenvalues $\lambda_j^\varepsilon$.

LEMMA 2.2. *There exists $C > 0$ such that for $1 \le j \le n$,*

(16)
$$\left|\lim_{\varepsilon \to 0} \varepsilon \ln \lambda_j^\varepsilon - [U(r_j, x_j) - U(r_j, -1)]\right| \le C\delta.$$

PROOF. Fix $1 \le j \le n$. Define the pseudopotential corresponding to the drift $R_j$ by

$$V_j(x, z) = \inf\left\{\frac{1}{2}\int_0^t \left(\phi_s' + \frac{R_j}{2}(\phi_s)\right)^2 ds, \phi_0 = x, \phi_t = z, t > 0\right\},$$

where $\phi$ stands for absolutely continuous functions defined on the time interval $[0, t]$. Since due to our assumptions $\underline{m}_j$ is the only local minimum of the potential corresponding to $R_j$ on $]-\infty, 0[$, the sharpened form of the exit time theorem of Freidlin–Wentzell (see [4], Theorem 1.1) implies that

$$\lim_{\varepsilon \to 0} \varepsilon \ln \lambda_j^\varepsilon = -V_j(\underline{m}_j, x_j).$$

Let us estimate the pseudopotential. We have

$$V_j(\underline{m}_j, x_j) = \int_{\underline{m}_j}^{x_j} R_j(\theta)\, d\theta$$

$$= \int_{-1}^{x_j} \frac{\partial U}{\partial x}(r_j, \theta)\, d\theta + \int_{\underline{m}_j}^{x_j} \left( R_j(\theta) - \frac{\partial U}{\partial x}(r_j, \theta) \right) d\theta$$

$$- \int_{[-1, \underline{m}_j]} \frac{\partial U}{\partial x}(r_j, \theta)\, d\theta.$$

Continuity of $\frac{\partial U}{\partial x}$ in $(t, x)$ entails the existence of $C_1 < 0$ such that

$$\left| \int_{[-1, \underline{m}_j]} \frac{\partial U}{\partial x}(r_j, \theta)\, d\theta \right| \leq C_1 \delta.$$

To estimate the second remainder term, recall that the mesh $\gamma$ was chosen to produce at most $\delta$ as modulus of continuity of $\frac{\partial U}{\partial x}$ [see (12)], and that $R_j$ is also at most a distance $\delta$ away [see (13)]. We therefore obtain

$$\left| \int_{\underline{m}_j}^{x_j} \left( R_j(\theta) - \frac{\partial U}{\partial x}(r_j, \theta) \right) d\theta \right| \leq 2\delta.$$

Hence

(17) $$\left| \lim_{\varepsilon \to 0} \varepsilon \ln \lambda_j^\varepsilon - [U(r_j, x_j) - U(r_j, -1)] \right| \leq (2 + C_1)\delta.$$

The asserted asymptotic result follows. $\square$

As a consequence we obtain an upper bound for the exponential convergence rate for the exit time from the domain of attraction of the potential well at $-1$.

PROPOSITION 2.1. *Let $x < 0$ and let $\Gamma$ be a compact subset of (10). Then there exists $0 < h_0 < \inf_{\mu \in \Gamma} a_\mu$ such that for $h \leq h_0$,*

(18) $$\limsup_{\varepsilon \to 0} \varepsilon \ln \mathbb{P}_x(\tau_1(X) \leq (a_\mu - h)T) \leq \mu - D_{-1}(a_\mu - h)$$

*uniformly for $\mu \in \Gamma$.*

PROOF. According to what has been proved, there are constants $\varepsilon_0 > 0$ and $K > 0$ such that for $\varepsilon \leq \varepsilon_0, \mu \in \Gamma$,

$$\mathbb{P}_x(\tau_1(X) \leq (a_\mu - h)T)$$

(19) $$\leq \sum_{j=1}^{n-1} \mathbb{P}_{x_{j-1}}(\tau_{x_j}(X^j) \leq \gamma T) + \mathbb{P}_{x_{n-1}}(\tau_1(X^n) \leq \gamma T)$$

$$\leq K\gamma T \sum_{j=1}^{n} \lambda_j^\varepsilon + n e^{-c/\varepsilon}.$$

Taking logarithms on both sides, multiplying by $\varepsilon$ and using the equation

$$\lim_{\varepsilon \to 0} \varepsilon \ln(f(\varepsilon) + g(\varepsilon)) = \max\left[\lim_{\varepsilon \to 0} \varepsilon \ln(f(\varepsilon)), \lim_{\varepsilon \to 0} \varepsilon \ln(g(\varepsilon))\right]$$

for two positive functions $f$ and $g$, we may apply Lemma 2.2 to get

$$\limsup_{\varepsilon \to 0} \varepsilon \ln \mathbb{P}_x(\tau_1(X) \leq (a_\mu - h)T)$$
$$\leq [\max[\mu - [U(r_j, x_j) - U(r_j, -1)] : 1 \leq j \leq n] + C\delta] \vee (-c)$$
$$\leq [\max[\mu - [U(r_j, x_1) - U(r_j, -1)] : 1 \leq j \leq n] + C\delta] \vee (-c).$$

Recalling the definition of $a_\mu$ and that $x_1 < 0$ is arbitrary, we may conclude

$$\limsup_{\varepsilon \to 0} \varepsilon \ln \mathbb{P}_x(\tau_1(X) \leq (a_\mu - h)T)$$
$$\leq [\max[\mu - D_{-1}(r_j) : 1 \leq j \leq n] + C\delta] \vee (-c)$$
$$\leq [(\mu - D_{-1}(a_\mu - h)) + C\delta] \vee (-c)$$

uniformly for $\mu \in \Gamma$. Now choose $h_0 > 0$ small enough so that for $h \leq h_0$ we have

$$\inf_{\mu \in \Gamma}(\mu - D_{-1}(a_\mu - h)) > -c.$$

Finally, since $\delta$ is arbitrary, we may let $\delta$ tend to zero. This way we obtain the desired upper bound for the exponential rate. □

2.2. *Lower bound for the exponential exit rate.* In the second step of our approach, we shall establish lower bounds for the exponential rates at which the diffusion exits from the basin of attraction of $-1$. Let us first prove an auxiliary result. It states that the probability of exiting the interval $[l, 0]$ via $l$ is exponentially small with exponential order increasing in $|l|$, due to hypercontractivity. Recall the constants $K_1$ and $K_2$ from (8).

LEMMA 2.3. *There exist positive constants $C$ and $\varepsilon_0$ such that for $\varepsilon \leq \varepsilon_0$ and $l < x \wedge -K_1$ and $\mu > 0$ we have*

$$(20) \qquad \mathbb{P}_x(\tau_l(X) \leq T) \leq \frac{C}{\varepsilon} \exp\left(\frac{2K_2(l - x \wedge (-K_1)) + \mu}{\varepsilon}\right).$$

PROOF. We give arguments for the case $x < -K_1$, the other case being easier. Recalling that by (8) the gradient of $U$ is bounded below by $-K_2$ on $\mathbf{R}_+ \times ]-\infty, -K_1]$, we may compare the diffusion $X$ with the diffusion $Z$ on the interval $]-\infty, x]$ reflected at $x$ with constant drift equal to $K_2$. It can be given by the SDE

$$dZ_t = K_2 \, dt + dL_t + \sqrt{2\varepsilon} \, dW_t,$$

where $Z_0 = x$ and $L$ is an increasing process satisfying $\int_0^t (Z_s - x)\,dL_s = 0$, $t \geq 0$. See, for example, [19]. The comparison clearly yields

$$\mathbb{P}_x(\tau_l(X) \leq T) \leq \mathbb{P}_x(\tau_l(Z) \leq T). \tag{21}$$

By Chebyshev's inequality,

$$\mathbb{P}_x(\tau_l(Z) \leq T) \leq e\mathbb{E}_x\left[\exp\left(-\frac{1}{T}\tau_l(Z)\right)\right]. \tag{22}$$

Let $\varphi(y) = \mathbb{E}_y[\exp(-\frac{1}{T}\tau_l(Z))]$, $l \leq y \leq x$. Our task consists in an estimation of $\varphi(x)$. According to the Feynman–Kac and Dynkin formulae $\varphi$ solves the boundary value problem

$$\varepsilon\varphi'' + K_2\varphi' - \frac{1}{T}\varphi = 0 \quad \text{on } ]l, x[,$$

$$\varphi'(x) = 0, \quad \varphi(l) = 1.$$

The eigenvalues of the differential equation are determined by the equation

$$\varepsilon\lambda^2 + K_2\lambda - \frac{1}{T} = 0,$$

hence by $\lambda^{\pm} = \frac{1}{2\varepsilon}[-K_2 \pm \sqrt{K_2^2 + \frac{4\varepsilon}{T}}]$. Taking the boundary conditions into account leads to the equation

$$\varphi(y) = \frac{\lambda^+ e^{\lambda^+ x + \lambda^- y} - \lambda^- e^{\lambda^- x + \lambda^+ y}}{\lambda^+ e^{\lambda^+ x + \lambda^- l} - \lambda^- e^{\lambda^- x + \lambda^+ l}}, \quad y \in [l, x].$$

Neglecting the second term in the denominator of the fraction, we obtain

$$\varphi(x) \leq \frac{\sqrt{K_2^2 + 4\varepsilon/T}}{\varepsilon\lambda^+ e^{\lambda^-(l-x)}}.$$

Now for $\varepsilon$ small enough, $\lambda^+ \geq \frac{1}{K_2 T}$ and $\lambda^- \leq \frac{-K_2}{\varepsilon}$. Therefore, for $\varepsilon$ small enough there exists a constant $C_0 > 0$ independent of $\varepsilon$ and $T$ such that

$$\varphi(x) \leq C_0 K_2^2 \frac{T}{\varepsilon} \exp\left[\frac{K_2(l-x)}{\varepsilon}\right] = \frac{C_0 K_2^2}{\varepsilon} \exp\left[\frac{K_2(l-x) + \mu}{\varepsilon}\right].$$

This implies the desired inequality (20). □

We shall continue to use the partition $(I_j : 1 \leq j \leq n)$ of the interval $[0, a_\mu - h]$ of the preceding section. This time, we shall compare with homogeneous diffusions by freezing the potential derivative at an upper level, which results in working with the drifts $V^{I_j}$ and the diffusions $\underline{X}^{I_j}$. Again, these drifts may fail to possess the regularity properties required to apply

Theorem 1.1. For this reason we may choose smoothed versions $V_j \in \mathcal{C}^1(\mathbb{R})$ of the potentials with corresponding diffusion processes $Y^j$ satisfying

(23)
$$\begin{aligned}
&V_j \geq V_{I_j} \geq V_j + \delta, \\
&\text{there are } \underline{m}_j \in \,]-1-\delta, -1+\delta[\,, s_j \in \,]-\delta, \delta[\,, \overline{m}_j \in \,]1-\delta, 1+\delta[ \\
&\text{such that } V_j'|_{]-\infty,\underline{m}_j[\cup]s_j,\overline{m}_j[} < 0, V_j'|_{]\underline{m}_j,s_j[\cup]\overline{m}_j,\infty[} > 0, \\
&V_j''(\underline{m}_j) > 0, \\
&\text{there are constants } K_1, K_2 > 0 \text{ such that } V_j(x) \leq -K_2 \text{ for } x \leq -K_1, \\
&V_j(x) \geq K_2 \text{ for } x \geq K_1.
\end{aligned}$$

To deduce a lower estimate, we shall compare $X$ via $X^{I_n}$ with $Y^n$ on the interval $I_n$ in the scale $T$. We may write for $l < x$

(24)
$$\begin{aligned}
\mathbb{P}_x&(\tau_1(X) \leq (a_\mu - h)T) \\
&\geq \mathbb{P}_x(\tau_1(X) \leq r_n T, \tau_1(X) \wedge \tau_l(X) \geq r_{n-1}T) \\
&\geq \mathbb{E}_x(\mathbb{1}_{\{\tau_1(X) \wedge \tau_l(X) \geq r_{n-1}T\}} P_{X_{r_{n-1}T}}(\tau_1(X_{\cdot + r_{n-1}T}) \leq \gamma T)) \\
&\geq \mathbb{P}_l(\tau_1(Y^n) \leq \gamma T) \times \mathbb{P}_x(\tau_1(X) \wedge \tau_l(X) \geq r_{n-1}T).
\end{aligned}$$

As a consequence of $\mu - D_{-1}((a_\mu - h)T) < 0$ and the arguments presented in Section 2.1, we note that uniformly on our compact set $\Gamma$

$$\lim_{\varepsilon \to 0} \mathbb{P}_x(\tau_1(X) \leq r_{n-1}T) = 0.$$

This clearly implies that there is $\varepsilon_0 > 0$ and a constant $C > \frac{1}{2}$ such that for $\varepsilon \leq \varepsilon_0$,

(25)
$$\mathbb{P}_x(\tau_1(X) \geq r_{n-1}T) \geq C.$$

Moreover, by Lemma 2.3, for $l$ small enough, there exists $\varepsilon_1 > 0$ such that for $\varepsilon \leq \varepsilon_1$,

$$\mathbb{P}_x(\tau_1(X) \geq r_{n-1}T) \geq C.$$

Hence for $\varepsilon$ small enough we have

$$\mathbb{P}_x(\tau_1(X) \wedge \tau_l(X) \geq r_{n-1}T) \geq C - \tfrac{1}{2} > 0.$$

It therefore remains to find lower bounds for $\mathbb{P}_l(\tau_1(Y^n) \leq \gamma T)$.

We may now apply the same arguments as those developed in Section 2.1. We just have to use Lemma 2.2 for the eigenvalues of the operator

$$\mathcal{L}^\varepsilon u = \varepsilon u'' - V_n u'$$

with Dirichlet boundary conditions at 0 in the sense of lower bounds uniformly on the compact $\Gamma$. As a consequence we obviously obtain, with a constant $C > 0$ independent of $\Gamma$,

$$\liminf_{\varepsilon \to 0} \varepsilon \ln \mathbb{P}_l(\tau_1(Y^n) \leq \gamma T) \geq (\mu - D_{-1}(a_\mu - h) - C\delta) \vee (-c),$$

uniformly for $\mu \in \Gamma$. Let us now choose $h_0 > 0$ and $\delta_0 > 0$ small enough such that for $h \leq h_0, \delta \leq \delta_0$ we have $\mu - D_{-1}((a_\mu - h)T) - C\delta \geq -c$. Since $\delta$ is arbitrary, we obtain

$$\liminf_{\varepsilon \to 0} \varepsilon \ln \mathbb{P}_l(\tau_1(Y^n) \leq \gamma T) \geq (\mu - D_{-1}((a_\mu - h)T)),$$

uniformly for $\mu \in \Gamma$. Recalling (24) and (25), we finally obtain

$$\liminf_{\varepsilon \to 0} \varepsilon \ln \mathbb{P}_x(\tau_1(X) \leq (a_\mu - h)T) \geq \mu - D_{-1}((a_\mu - h)T),$$

uniformly on $\Gamma$. With this result we have established the desired lower bound for the exponential exit rate.

PROPOSITION 2.2. *Let $x < 0$, and let $\Gamma$ be a compact subset of* (10). *Then there exists $0 < h_0 < \inf_{\mu \in \Gamma} a_\mu$ such that for $h \leq h_0$,*

$$(26) \qquad \lim_{\varepsilon \to 0} \varepsilon \ln \mathbb{P}_x(\tau_1(X) \leq (a_\mu - h)T) \geq \mu - D_{-1}(a_\mu - h)$$

*uniformly for $\mu \in \Gamma$.*

2.3. *The exponential smallness of the rate of too long transitions.* Having proved (11) in the preceding two propositions, the second aim of this section is to show that the exponential rate at which $\tau_1(X)$ exceeds $(a_\mu + h)T$ is arbitrarily small. In fact, we shall make precise that the rate at which transitions happen which take at least as long as $(a_\mu + h)T$ vanishes to all exponential orders, for $h > 0$ arbitrary.

PROPOSITION 2.3. *Let $x < 0$, and let $\Gamma$ be a compact subset of* (10) *not containing $D_{-1}(0)$. Then there exists $h_0 > 0$ such that for all $0 < h \leq h_0$ and $\mu \in \Gamma$ we have*

$$\limsup_{\varepsilon \to 0} \varepsilon \ln \mathbb{P}_x(\tau_1(X) \geq (a_\mu + h)T) = -\infty.$$

PROOF. Let $\delta > 0$ and $h > 0$ be given. Let $\Gamma$ be a compact subset of (10). First, for $l < x$ we may write

$$(27) \qquad \mathbb{P}_x(\tau_1(X) \geq (a_\mu + h)T) \leq \mathbb{P}_x(\tau_1(X) \wedge \tau_l(X) \geq (a_\mu + h)T)$$
$$(28) \qquad \qquad \qquad + \mathbb{P}_x(\tau_l(X) \leq (a_\mu + h)T).$$

To estimate the second term on the right-hand side of (27), we employ Lemma 2.3. In fact, for $l < x \wedge (-K_1)$ we have

$$\limsup_{\varepsilon \to 0} \varepsilon \ln \mathbb{P}_x(\tau_l(X) \leq (a_\mu + h)T) \leq 2K_2(l - x \wedge (-K_1)) + \sup_{\mu \in \Gamma} \mu.$$

Therefore

$$\lim_{l \to -\infty} \left[ \limsup_{\varepsilon \to 0} \varepsilon \ln \mathbb{P}_x(\tau_l(X) \leq (a_\mu + h)T) \right] = -\infty.$$

It therefore remains to estimate the first term on the right-hand side of (27) for $l$ small but fixed. Let $0 = r_0 < r_1 < \cdots < r_n = a_\mu + h$ be an equidistant partition of the interval $[0, a_\mu + h]$ of mesh $\gamma < \frac{h}{2}$ and denote $I_j = [r_{j-1}, r_j], 1 \leq j \leq n$. Then we have

$$\mathbb{P}_x(\tau_1(X) \wedge \tau_l(X) \geq (a_\mu + h)T)$$

$$= \mathbb{P}_x(\tau_0(X) \wedge \tau_l(X) \geq r_{n-1}T, \tau_1(X_{\cdot + r_{n-1}T}) \wedge \tau_l(X_{\cdot + r_{n-1}T}) \geq \gamma T)$$

$$\leq \mathbb{E}_x(\mathbb{1}_{\{\tau_0(X) \wedge \tau_l(X) \geq r_{n-1}T\}} \mathbb{P}_{X_{r_{n-1}T}}(\tau_1(\underline{X}^{I_n}) \geq \gamma T))$$

(29)
$$\leq \max_{y \in ]l, 0[} \mathbb{P}_y(\tau_1(\underline{X}^{I_n}) \geq \gamma T)$$

$$= \mathbb{P}_l(\tau_1(\underline{X}^{I_n}) \geq \gamma T)$$

$$\leq \mathbb{P}_l(\tau_1(X^n) \geq \gamma T).$$

Here, we compare the inhomogeneous diffusion $X$ on $I_n$ with the time-homogeneous one $X^{I_n}$ corresponding to drift $R_{I_n}$ and finally with $X^n$ subject to drift $R_n$ to be described below, and we use monotonicity of

$$y \mapsto \mathbb{P}_y(\tau_1(\underline{X}^{I_n}) \geq \gamma T).$$

We assume $\gamma$ to be small enough to ensure

(30) $$\sup_{s,t \in [r_{n-1}, r_n]} \sup_{x \in [l,0]} \left| \frac{\partial U}{\partial x}(t, x) - \frac{\partial U}{\partial x}(s, x) \right| \leq \delta.$$

We may choose the drift $R_n$ to satisfy

(31) $R_n \leq R_{I_n} \leq R_n + \delta$,
there are $\underline{m}_n \in ]-1-\delta, -1+\delta[, s_n \in ]-\delta, \delta[, \overline{m}_n \in ]1-\delta, 1+\delta[$
such that $R_n|_{]-\infty, \underline{m}_n[ \cup ]s_n, \overline{m}_n[} < 0, R_n|_{]\underline{m}_n, s_n[ \cup ]\overline{m}_n, \infty[} > 0$,
$R'_n(\underline{m}_n) > 0$,
there are constants $K_1, K_2 > 0$ such that $R_n(x) \leq -K_2$ for $x \leq -K_1$,
$R_n(x) \geq K_2$ for $x \geq K_1$.

We are ready to apply Theorem 1.1, this time for $Q = R_n, d = -\delta$ to obtain that

(32) $$\mathbb{P}_l(\tau_1(X^n) \geq \gamma T) \leq e^{-\lambda_n^\varepsilon \gamma T}(1 + e^{-c/\varepsilon})$$

uniformly on $\Gamma$. Here $\lambda_n^\varepsilon$ stands for the principal eigenvalue of the operator $\mathcal{L}_n^\varepsilon$ defined by

$$\mathcal{L}_n^\varepsilon u = \varepsilon u'' - R_n u'$$

with Dirichlet boundary conditions at $-\delta$. The asymptotic properties of $\lambda_n^\varepsilon$ can be deduced in a similar way to Lemma 2.2. We estimate the pseudopotential corresponding to $R_n$ taking (30) into account. We obtain that there

exists $C > 0$ such that

(33) $$\left|\lim_{\varepsilon \to 0} \varepsilon \lambda_n^\varepsilon T - (\mu - D_{-1}(r_n))\right| \leq C\delta$$

uniformly on $\Gamma$. Now recall that due to the choice of $\gamma$, we have $a_\mu + \frac{h}{2} < r_n$, hence

$$\mu - D_{-1}(r_n) > 0, \qquad \mu \in \Gamma,$$

and by compactness of $\Gamma$ even,

$$\inf_{\mu \in \Gamma}[\mu - D_{-1}(r_n)] > 0.$$

This in turn implies that

$$\lim_{\varepsilon \to 0} \lambda_n^\varepsilon T = \infty$$

uniformly on $\Gamma$. But due to (32) we are allowed to conclude

$$\lim_{\varepsilon \to 0} \mathbb{P}_l(\tau_1(X^n) \geq \gamma T) = -\infty$$

uniformly in $\Gamma$. According to (29), this completes the proof. $\square$

Finally, we may summarize the results of Sections 2.1 and 2.2, and state the main result on asymptotic exponential decay rates of transition probabilities.

THEOREM 2.1. *Let $x < 0$, and let $\Gamma$ be a compact subset of* (10). *For $\varepsilon > 0$, $\mu \geq 0$ let $T = \exp(\frac{\mu}{\varepsilon})$. Then there exists $h_0 > 0$ such that for $h \leq h_0$,*

(34) $$\lim_{\varepsilon \to 0} \varepsilon \ln \mathbb{P}_x(\tau_1(X) \notin [(a_\mu - h)T, (a_\mu + h)T]) = \mu - D_{-1}(a_\mu - h)$$

*uniformly for $\mu \in \Gamma$.*

PROOF. This is an immediate consequence of Propositions 2.1–2.3. $\square$

**3. Stochastic resonance in a double-well potential.** Let us now turn to the main subject of this paper, a characterization of the notion of *stochastic resonance*. Let us recall that we look for a characterization of the concept of *optimal periodic tuning* which is extensively studied in the physics literature by notions such as the signal-to-noise ratio or the spectral power amplification (see [18]). Let us also remark that this concept implicitly uses and refines the concept of stochastic resonance studied by Freidlin [7] which paraphrases the ability of periodically perturbed stochastic systems to follow the periodic excitation in the small noise limit, and exhibit *quasi-periodic motion*. In more mathematical terms and the notation introduced before, we aim at choosing the noise intensity parameter $\varepsilon$ such that in the *large period limit* $T \to \infty$ the diffusion trajectories follow the periodic excitation of the system hidden in $U$ in an optimal way to be made precise. In Section 3.1 we shall show that a quality measure of goodness of periodic tuning is given by the exponential rate at which the first transition to the other well happens within a fixed interval around $a_\mu T$. In Section 3.2 we establish robustness of this notion of quality: we show that in the small noise limit the diffusion and its reduced model, a Markov chain living on a two-point state space, have the same resonance pattern.

3.1. *Transition probabilities as a measure of quality.* The local extrema of the depth functions $D_{\pm 1}$ of $U$ are supposed to be global, and $D_{\pm 1}$ is strictly increasing between its extrema. Recall that we work with exponential time scales $\mu$ related to the natural time $T$ by the equation $T = \exp(\frac{\mu}{\varepsilon})$. In this section, we have to work with scale functions depending on the starting well and eventually on arbitrary starting times. So we let

$$a_\mu^i(s) = \inf\{t \geq s : D_i(t) \leq \mu\}, \qquad i = \pm 1, \mu \geq 0.$$

The relevant time scales $\mu$ will be chosen from the intervals

$$I_i = \left]\inf_{t \geq 0} D_i(t), \sup_{t \geq 0} D_i(t)\right[, \qquad i = \pm 1.$$

Our aim is to observe periodic behavior of the diffusion. This will in principle mean that the process can travel from one well to the other and back on the time scales in which we let the diffusion run, but not instantaneously. So, on the one hand, we have to work on time scales on which it *does not get stuck in one of the wells* of the potential. On the other hand, the time scales we are concentrating on should also *not allow for chaotic behavior*, that is, immediate re-bouncing after changing the well.

To make these conditions mathematically precise, recall that transitions become possible as soon as the potential barrier $D_{\pm 1}$ becomes smaller than the time scale parameter $\mu$. Hence if $\mu > \inf_{t \geq 0} D_i(t)$, there is a time range during which the diffusion can leave the well centered at $i$. To not get stuck in one of them, the diffusion has to be able to leave both. This is guaranteed if

$$(35) \qquad \mu > \max_{i = \pm 1} \inf_{t \geq 0} D_i(t).$$

To avoid immediate re-bouncing, we have to assure that the diffusion cannot leave the domain of attraction of $-i$ at the moment it reaches it, coming from $i$. Suppose we consider the dynamics after time $s \geq 0$, and the diffusion is near $i$ at that time. Its first transition to the well at $-i$ occurs at time $a_\mu^i(s)T$, and it stays there for at least a little while if $D_{-i}(a_\mu^i(s))$ is bigger than $\mu$. This is equivalent to stating that for all $s \geq 0$ there exists $\delta > 0$ such that on $[a_\mu^i(s), a_\mu^i(s) + \delta]$ we have $\mu < D_{-i}$. But for $t$ shortly after $a_\mu^i(s)$, we always have $D_i(t) \leq \mu$ by the very definition of $a_\mu^i$. Hence our condition becomes equivalent to the following: for all $s \geq 0$ there exists $\delta > 0$ such that on $[a_\mu^i(s), a_\mu^i(s) + \delta]$ we have $\mu < \max_{i=\pm 1} D_i$. This in turn is more elegantly expressed by

$$(36) \qquad \mu < \inf_{t \geq 0} \max_{i=\pm 1} D_i(t).$$

See Figure 8.

We may summarize our search for an appropriate set of scale parameters $\mu$ for which periodicity in the diffusion behavior will occur. We call this set the "resonance interval" to indicate that we have to look for the scale of optimal periodicity, the *resonance scale*, in this interval. See [11] for the definition of the corresponding interval in the case of two-state Markov chains. The interval

$$I_R = \left] \max_{i=\pm 1} \inf_{t \geq 0} D_i(t), \inf_{t \geq 0} \max_{i=\pm 1} D_i(t) \right[$$

is called the *resonance interval* (see Figure 3). Let us pause for a moment at this point to compare our approach with Freidlin's [7] understanding of stochastic resonance by quasi-deterministic motion. In Freidlin's terms, stochastic resonance is given if the parameter $\mu$ exceeds the lower boundary of our resonance interval. Our concept of resonance stipulates to look for an optimal $\mu$ in the resonance interval at which in a sense to be made precise the quality of periodic tuning is optimal.

Let us now come to the discussion of the quality of periodic response of the stochastic system given by the diffusion, in dependence on the noise parameter $\varepsilon$ and the time scale parameter $\mu$ which according to the remarks made above has to be chosen in the resonance interval. To simplify things a little, let us assume that the depth functions are related by a phase $\phi \in \,]0,1[$, that is,

$$D_{-1}(t) = D_1(t + \phi), \qquad t \geq 0.$$

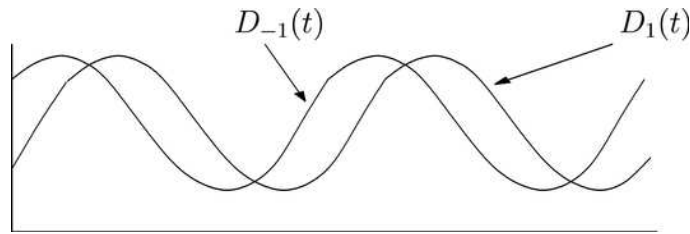

FIG. 8. *Depth functions in phase.*

Moreover, we assume that the diffusion starts in $-1$. There are many ways to describe optimality of periodic tuning. Imkeller and Pavlyukevich [12, 13] consider different measures of quality such as the *spectral power amplification*, the *energy*, the *energy-to-noise ratio*, the *out-of-phase measure*, the *entropy* and the *relative entropy*. The detailed study of the physicists' favorite measure, the spectral power amplification, based on the energy of the spectral component of the mean trajectory in equilibrium corresponding to the forcing frequency $\frac{2\pi}{T}$, shows one surprising defect: it is not robust as one passes from the diffusion to a reduced model described by a two-state Markov chain jumping with rates corresponding to the transition rates between the metasable states $\pm 1$ of the diffusion given by the potential minima. In fact, while the Markov chain's spectral power coefficient shows a pronounced peak for $\mu$ near an average well depth, the overwhelming influence of the diffusion's fluctuations in small neighborhoods of the potential wells, discovering very subtle details of the potential's geometry there, destroys this picture completely. Here we propose a notion of quality of periodic tuning which is based on the pure transition mechanism of the system between the domains of attraction of the double-well potential. Generalizing the approach of a study of optimal tuning for two-state Markov chain models (see [11]), we measure the quality of tuning by computing for varying time scale parameters $\mu$ the probability that, starting in $i$, the diffusion is transferred to $-i$ within the time window $[(a_\mu^i - h)T, (a_\mu^i + h)T]$ of width $2hT$. To find the *stochastic resonance point* for large $T$ we have to maximize this measure of quality in $\mu \in I_R$. The probability for transition within this window will be computed by the estimates of the preceding section. Uniformity of convergence to the exponential rates will enable us to maximize in $\mu$. Let us now make these ideas precise.

To make sure that the transition window makes sense at least for small $h$, we have to suppose that $a_\mu^i > 0, i = \pm 1$ for $\mu \in I_R$. This will be guaranteed if

$$(37) \qquad D_i(0) > \inf_{t \geq 0} \max_{i = \pm 1} D_i(t), \qquad i = \pm 1.$$

If this is not granted from the beginning, it suffices to start the diffusion a little later, in order to be sure that (37) is satisfied. Under (37), we call

$$(38) \quad M(\varepsilon, \mu) = \min_{i = \pm 1} \mathbb{P}_i(\tau_{-i}(X) \in [(a_\mu^i - h)T, (a_\mu^i + h)T]), \qquad \varepsilon > 0, \mu \in I_R,$$

*the transition probability for a time window of width $h$.*

We are prepared to state our main resonance result.

THEOREM 3.1.  *Let $\Gamma$ be a compact subset of $I_R$, and let $h_0 > 0$ be given according to Theorem* 2.1. *Then*

$$(39) \qquad \lim_{\varepsilon \to 0} \varepsilon \ln(1 - M(\varepsilon, \mu)) = \max_{i = \pm 1}\{\mu - D_i(a_\mu^i - h)\}$$

*uniformly for $\mu \in \Gamma$.*

Proof. This proposition is an obvious consequence of Propositions 2.1–2.3.

□

It is clear that for $h$ small the eventually existing global minimizer $\mu_R(h)$ of

$$I_R \ni \mu \mapsto \max_{i=\pm 1}\{\mu - D_i(a^i_\mu - h)\}$$

is a good candidate for our resonance point. But it still depends on $h$. To get rid of this dependence, we shall consider the limit of $\mu_R(h)$ as $h \to 0$.

Definition 3.1. Suppose that

$$I_R \ni \mu \mapsto \max_{i=\pm 1}\{\mu - D_i(a^i_\mu - h)\}$$

possesses a global minimum $\mu_R(h)$. Suppose further that

$$\mu_R = \lim_{h \to 0} \mu_R(h)$$

exists in $I_R$. We call $\mu_R$ the *stochastic resonance point* of the diffusion $X$ with time-periodic potential $U$.

We shall now show that the stochastic resonance point exists if one of the depth functions, and thus both, due to the phase lag, has a unique point of maximal decrease on the interval where it is strictly decreasing. See Figure 9.

Theorem 3.2. *Suppose that $D_1$ is twice continuously differentiable and has its global maximum at $t_1$, and its global minimum at $t_2$, where $t_1 < t_2$. Suppose further that there is a unique point $t_1 < s < t_2$ such that $D_1|_{]t_1,s[}$ is strictly concave, and $D_1|_{]s,t_2[}$ is strictly convex. Then $\mu_R = D_1(s)$ is the stochastic resonance point.*

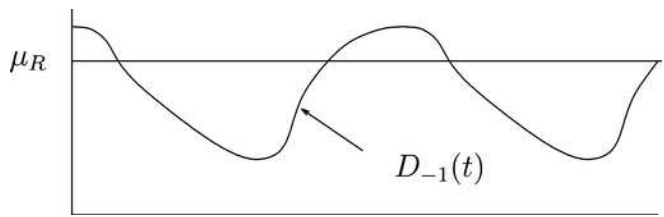

Fig. 9. *Point of maximal decrease.*

PROOF. First of all, note that there is $\psi \in ]0,1[$ such that $D_1 = D_{-1}(\cdot + \psi)$. As a consequence of this,

$$\max_{i=\pm 1}\{\mu - D_i(a_\mu^i - h)\} = \{\mu - D_1(a_\mu^1 - h)\}.$$

Write $a_\mu = a_\mu^1$ and recall that on the interval of decrease of $D_1$, $a_\mu = D_1^{-1}(\mu)$. Therefore, the differentiability assumption yields

$$1 = D_1'(a_\mu - h) \cdot a_\mu' = D_1'(a_\mu - h) \cdot \frac{1}{D_1'(a_\mu)}.$$

Our hypothesis concerning convexity and concavity of $D_1$ essentially means that $D_1''(s) = 0$, and $D_1''|_{]t_1,s[} < 0, D_1''|_{]s,t_2[} > 0$; in other words, that $\mu \mapsto D_1'(a_\mu)$ has a local maximum at $a_\mu = s$. Hence for $h$ small there exists a unique point $a_\mu(h)$ such that

$$D_1'(a_\mu(h) - h) = D_1'(a_\mu(h))$$

and

$$\lim_{h \to 0} a_\mu(h) = s.$$

To show that $a_\mu(h)$ corresponds to a minimum of the function

$$\mu \mapsto [\mu - D_1(a_\mu - h)],$$

we take the second derivative of this function at $a_\mu(h)$, which is given by

$$\frac{D_1'(a_\mu(h) - h)D_1''(a_\mu(h)) - D_1''(a_\mu(h) - h)D_1'(a_\mu(h))}{D_1'(a_\mu(h))}.$$

But $D_1'(a_\mu(h)), D_1'(a_\mu(h) - h) < 0$, whereas $D_1''(a_\mu(h) - h) > 0, D_1''(a_\mu(h)) < 0$. This clearly implies that $a_\mu(h)$ corresponds to a minimum of the function. But by definition, as $h \to 0$, $a_\mu(h) \to s$. Therefore, finally, $D_1(s)$ is the stochastic resonance point. □

To illustrate our results, we next discuss an example.

EXAMPLE. Let us consider the double-well potential

$$U(t,x) = \frac{x^6}{6} - \cos\left\{2\pi\left(t - \frac{1}{4} + \psi\mathrm{sgn}(x)\right)\right\}\left(\frac{x^5}{5} - \frac{x^3}{3}\right) - \frac{x^2}{2},$$

with $T = \exp(\frac{\mu}{\varepsilon})$ and $\psi \in [0, \frac{1}{4}[$. See Figures 10 and 11.

$U$ satisfies all the assumptions required for potentials above, in particular

$$\frac{\partial U}{\partial x}(t,x) = 0 \quad \text{iff } x \in \{-1, 0, 1\}.$$

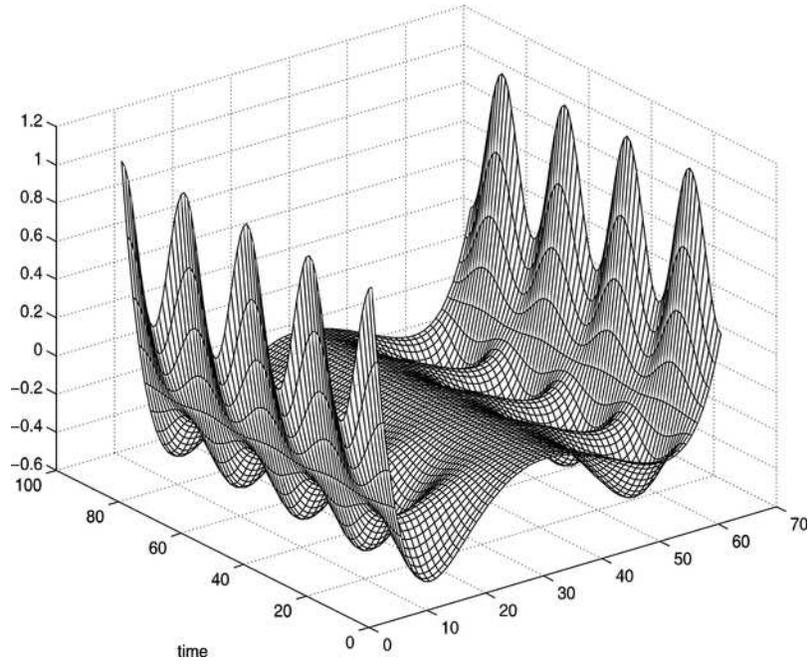

Fig. 10. *Double-well potential (case: $\psi = 0$).*

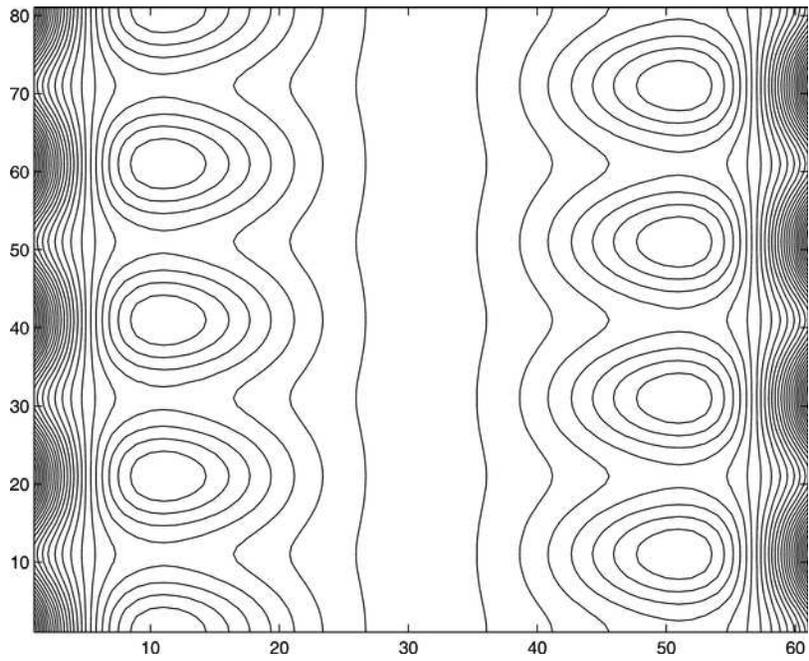

Fig. 11. *Level sets of the potential.*

$-1$ and $1$ are the metasable states of the potential, and $0$ is the saddle point. We can then compute the barrier height of the two wells. For $i \in \{-1, 1\}$,

$$D_i(t) = \frac{2}{3} - i\frac{4}{15}\cos\left(2\pi t + i2\pi\psi - \frac{\pi}{2}\right).$$

Let us note that $D_1(t) = D_{-1}(t + 2\psi + \frac{1}{2})$. Since $\psi \in [0, \frac{1}{4}[$, we are in the phase case with $\phi = 2\psi + 1/2$. The resonance interval is then given by

$$I_R = \left]\frac{2}{5}, \frac{2}{3} - \frac{4}{15}\cos\left(\frac{\pi}{2} - 2\pi\psi\right)\right[.$$

In the symmetric case, that is, if $\psi = 0$, we obtain $I_R = ]2/5, 2/3[$. Let us now compute the optimal tuning scale applying Theorem 2.1. We obtain

$$a_\mu^{-1} = \frac{1}{2\pi}\arccos\left(\frac{15}{2}\left(\frac{\mu}{2} - \frac{1}{3}\right)\right) + \psi + \frac{1}{4}.$$

Hence, for $h > 0$ small enough,

$$F(\mu) = \mu - D_{-1}(a_\mu^{-1} - h)$$

$$= \mu - \frac{2}{3} - \frac{4}{15}\cos\left(\arccos\left\{\frac{15}{2}\left(\frac{\mu}{2} - \frac{1}{3}\right)\right\} - 2\pi h\right)$$

$$= \left(\mu - \frac{2}{3}\right)(1 - \cos 2\pi h) - \frac{4}{15}\sin(2\pi h)\sqrt{1 - \left(\frac{15}{2}\left(\frac{\mu}{2} - \frac{1}{3}\right)\right)^2}.$$

Let us recall that $F$ does not depend on the phase which implies

$$\mu - D_{-1}((a_\mu^{-1} - h)T) = \mu - D_1((a_\mu^1 - h)T).$$

Hence, to obtain optimal tuning, it suffices to compute the minimum of $F$ for $\mu \in I_R$. Differentiating $F$, we obtain

$$F'(\mu) = 1 - \cos(2\pi h) + \frac{15}{2}\sin(2\pi h)\frac{\mu/2 - 1/3}{\sqrt{1 - (15\mu/4 - 15/6)^2}}.$$

Hence $F$ attains its minimum for

$$\mu_R(h) = \frac{2}{3} - \frac{2\sqrt{2}}{15}\sqrt{1 - \cos(2\pi h)}$$

and

$$\mu_R = \lim_{h \to 0} \mu_R(h) = \frac{2}{3}.$$

Thus we obtain that $\mu_R$ is the stochastic resonance point if $\mu_R \in I_R$, that is, if the phase is near to $\frac{1}{2}$, that is, if $\psi$ is close to $0$. In the other case, the optimal tuning rate on every interval $[a, b] \subset I_R$ is given by the upper bound $b$.

3.2. *The robustness of stochastic resonance based on transition windows.*
In the small noise limit $\varepsilon \to 0$, it seems reasonable to assume that the periodicity properties of the diffusion trajectories caused by the periodic forcing due to the potential, are essentially captured by a simpler, reduced stochastic process: a continuous-time Markov chain which just jumps between two states representing the bottoms of the wells of the double-well potential at rates corresponding to the transition mechanism of the diffusion. This is just the reduction idea ubiquitous in the physics literature, and explained, for example, in [15]. In [13] it is found that this idea may conflict with the intrawell fluctuations of the diffusion if the quality of periodic tuning is measured by concepts using spectral decompositions of the trajectories. We shall now show that in the small noise limit both models, diffusion and Markov chain, produce the same resonance picture, if quality of periodic tuning is measured by transition rates as discussed in Section 3.1.

We first have to describe the reduced model. Let $U$ be a time-dependent potential function generating the potential diffusions of the preceding section. Recall that the depth functions of the potential minima satisfy $D_1(t) = D_{-1}(t+\phi)$, $t \geq 0$, with phase shift $\phi \in ]0,1[$. So, let us consider a time-continuous Markov chain $\{Y_t, t \geq 0\}$ taking values in the state space $\{-1,1\}$ with initial data $Y_0 = -1$. Suppose the infinitesimal generator is given by

$$G = \begin{pmatrix} -\varphi\left(\frac{t}{T}\right) & \varphi\left(\frac{t}{T}\right) \\ \psi\left(\frac{t}{T}\right) & -\psi\left(\frac{t}{T}\right) \end{pmatrix},$$

where $\psi(t) = \varphi(t+\phi)$, $t \geq 0$, and $\varphi$ is a 1-periodic function describing a rate which just produces the transition dynamics of the diffusion between the potential minima $\pm 1$, that is,

(40) $$\varphi(t) = \exp\left(-\frac{D_{-1}(t)}{\varepsilon}\right), \qquad t \geq 0.$$

Note that by choice of $\phi$,

(41) $$\psi(t) = \exp\left(-\frac{D_1(t)}{\varepsilon}\right), \qquad t \geq 0.$$

Transition probabilities for the Markov chain thus defined are easily computed. See ([11], Section 2). For example, the probability density of the first transition time $\sigma_i(Y)$ is given by

$$p(t) = \varphi(t) \exp\left(-\int_0^t \varphi(s)\,ds\right) \qquad \text{if } i = -1,$$

$$q(t) = \varphi(t+\phi) \exp\left(-\int_0^t \varphi(s+\phi)\,ds\right) \qquad \text{if } i = 1,$$

$t \geq 0$. Equation (42) can be used to obtain results on exponential rates of the transition times $\sigma_i(Y)$ if starting from $-i$, $i = \pm 1$. We summarize them and apply them to the measure of quality of periodic tuning in case (37):

(42) $\quad N(\varepsilon, \mu) = \min_{i=\pm 1} \mathbb{P}_i(\sigma_{-i}(Y) \in [(a_\mu^i - h)T, (a_\mu^i + h)T]), \qquad \varepsilon > 0, \mu \in I_R,$

which is called *transition probability for a time window of width h* for the Markov chain.

Here is the asymptotic result obtained from a slight modification of Theorems 3 and 4 of [11] which consists of allowing more general depth functions than the sinusoidal ones used there and requires just the same proof.

THEOREM 3.3. *Let $\Gamma$ be a compact subset of $I_R$ and let $h_0 < \sup(a_\mu^{-1}, T/2 - a_\mu^{-1})$. Then for $0 < h \leq h_0$,*

(43) $\qquad \lim_{\varepsilon \to 0} \varepsilon \ln(1 - N(\varepsilon, \mu)) = \max_{i=\pm 1} \{\mu - D_i(a_\mu^i - h)\}$

*uniformly for $\mu \in \Gamma$.*

It is clear from Theorem 3.3 that the reduced Markov chain $Y$ and the diffusion process $X$ have exactly the same resonance behavior. Of course, we may define the *stochastic resonance point* for $Y$ just as we did for $X$. So the following final robustness result holds true.

THEOREM 3.4. *The resonance points of $X$ with periodic potential $U$ and of $Y$ with exponential transition rate functions $D_{\pm 1}$ coincide.*

## REFERENCES

[1] BENZI, R., PARISI, G., SUTERA, A. and VULPIANI, A. (1983). A theory of stochastic resonance in climatic change. *SIAM J. Appl. Math.* **43** 565–578. MR700532
[2] BERGLUND, N. and GENTZ, B. (2002). Metastability in simple climate models: Pathwise analysis of slowly driven Langevin equations. *Stoch. Dyn.* **2** 327–358. MR1943556
[3] BOVIER, A., ECKHOFF, M., GAYRARD, V. and KLEIN, M. (2002). Metastability in reversible diffusion processes I. Sharp asymptotics for capacities and exit times. Unpublished manuscript. MR1911735
[4] BOVIER, A., GAYRARD, V. and KLEIN, M. (2002). Metastability in reversible diffusion processes II. Precise asymptotics for small eigenvalues. Unpublished manuscript.
[5] DAY, M. V. (1983). On the exponential exit law in the small parameter exit problem. *Stochastics* **8** 297–323. MR693886
[6] DAY, M. V. (1987). Recent progress on the small parameter exit problem. *Stochastics* **20** 121–150. MR877726
[7] FREIDLIN, M. (2000). Quasi-deterministic approximation, metastability and stochastic resonance. *Phys. D* **137** 333–352. MR1740101

[8] FREIDLIN, M. and WENTZELL, A. (1998). *Random Perturbations of Dynamical Systems*, 2nd ed. Springer, New York. MR1652127
[9] FREUND, J. A., NEIMAN, A. B. and SCHIMANSKY-GEIER, L. (2001). Stochastic resonance and noise-induced phase coherence. In *Stochastic Climate Models* (P. Imkeller and von J.-S. Storch, eds.) 309–323. Birkhäuser, Basel. MR1948303
[10] GAMMAITONI, L., HÄNGGI, P., JUNG, P. and MARCHESONI, F. (1998). Stochastic resonance. *Rev. Mod. Phys.* **70** 223–287.
[11] HERRMANN, S. and IMKELLER, P. (2002). Barrier crossings characterize stochastic resonance. *Stoch. Dyn.* **2** 413–436. MR1943561
[12] IMKELLER, P. and PAVLYUKEVICH, I. (2001). Stochastic resonance in two-state Markov chains. *Arch. Math. (Basel)* **77** 107–115. MR1845680
[13] IMKELLER, P. and PAVLYUKEVICH, I. (2002). Model reduction and stochastic resonance. *Stoch. Dyn.* **2** 463–506. MR1949298
[14] KRAMERS, H. A. (1940). Brownian motion in a field of force and the diffusion model of chemical reactions. *Phys. VII* **4** 284–304. MR2962
[15] MCNAMARA, B. and WIESENFELD, K. (1989). Theory of stochastic resonance. *Phys. Rev. A* **39** 4854–4869.
[16] MILSTEIN, G. M. and TRETJAKOV, M. V. (2000). Numerical analysis of noise-induced regular oscillations. *Phys. D* **140** 244–256. MR1757090
[17] NICOLIS, C. (1982). Stochastic aspects of climatic transitions—responses to periodic forcing. *Tellus* **34** 1–9. MR646774
[18] PAVLYUKEVICH, I. (2002). Stochastic resonance. Ph.D. thesis, Humboldt-Univ. zu Berlin.
[19] REVUZ, D. and YOR, M. (1999). *Continuous Martingales and Brownian Motion*, 3rd ed. Springer, Berlin. MR1725357

INSTITUT DE MATHÉMATIQUES ELIE CARTAN
UNIVERSITÉ HENRI POINCARÉ NANCY I
B.P. 239
54506 VANDOEUVRE-LÈS-NANCY CEDEX
FRANCE

INSTITUT FÜR MATHEMATIK
HUMBOLDT-UNIVERSITÄT ZU BERLIN
UNTER DEN LINDEN 6
10099 BERLIN
GERMANY
E-MAIL: imkeller@mathematik.hu-berlin.de